\documentclass[11pt,letterpaper]{article}

\usepackage{natbib, graphicx, amssymb, bm}

\def\a{\alpha}
\def\b{\beta}

\def\d{\delta}
\def\e{\epsilon}
\def\f{\phi} \def\vf{\varphi}
\def\g{\gamma}

\def\l{\lambda}
\def\m{\mu}
\def\n{\nu}

\def\q{\theta}
\def\r{\rho}
\def\s{\sigma}
\def\t{\tau}

\def\G{\Gamma}

\def\O{\Omega}

\def\la{\left}
\def\ra{\right}
\def\pa{\partial}
\def\del{\nabla}
\def\inf{\infty}

\def\abs#1{\left| #1\right|}

\def\Hat#1{\widehat{#1}}
\def\Tilde#1{\widetilde{#1}}

\newcommand{\ms}{\scriptscriptstyle}

\def\ba{\begin{array}}
\def\ea{\end{array}}
\def\be{\begin{equation}}
\def\ee{\end{equation}}
\def\bdm{\begin{displaymath}}
\def\edm{\end{displaymath}}
\def\bea{\begin{eqnarray}}
\def\eea{\end{eqnarray}}

\def\lb{\label}
\def\sp{~~~}

\begin{document}

\title{Heat Equation on the Cone and\\the Spectrum of the Spherical Laplacian}
\author{B S Balakrishna\footnote{balakbs2@gmail.com}}
\date{March 25, 2013}

\maketitle

\begin{abstract}

Spectrum of the Laplacian on spherical domains is analyzed from the point of view of the heat equation on the cone. The series solution to the heat equation on the cone is known to lead to a study of the Laplacian eigenvalue problem on domains on the sphere in higher dimensions. It is found that the solution leads naturally to a spectral function, a `generating function' for the eigenvalues and multiplicities of the Laplacian, expressible in closed form for certain domains on the sphere. Analytical properties of the spectral function suggest a simple scaling procedure for estimating the eigenvalues. Comparison of the first eigenvalue estimate with the available theoretical and numerical results for some specific domains shows remarkable agreement.

\end{abstract}

The heat equation on the $n$-dimensional cone has been addressed by various authors in the past and series solutions have been obtained. The $n=2$ solution was obtained by \citet{sxx}. The $n=3$ case was considered within the context of circular cones by \citet{cjx}. For higher dimensions, the applicable solution has been presented by \citet{cxx}. The radial component of the series solution involves modified Bessel functions while the angular component involves the eigenvalues and the eigenfunctions of the Laplacian on a domain on the $n-1$ dimensional sphere. It is found that the solution leads naturally to a spectral function, a `generating function' for the eigenvalues and their multiplicities, expressible in closed form for certain domains on the sphere such as the octant triangle on the two-sphere and analogous ones on higher dimensional spheres. Analytical properties of the spectral function suggest a simple scaling procedure to estimate the first few eigenvalues for related domains. The estimates for some specific domains are found to be in excellent agreement with the available theoretical and numerical results.

\section{Heat Equation on the Cone}
\lb{Hec}

Consider ${\mathbb R}^n$ with coordinate vector $\bm{x}$, conveniently split into radial and angular parts, $r$ and $\bm{\hat{r}}$ (when $r\ne 0$), as
\be
r^2 = \sum_i x_i^2, \sp \bm{\hat{r}} = \frac{\bm{x}}{r}.
\lb{rxy}
\ee
Unit radial vectors $\bm{\hat{r}}$ trace out the $n-1$ dimensional sphere $S^{n-1}$ at $r^2=1$. Consider a domain $\O^{n-1}$ of dimension $n-1$ on $S^{n-1}$ with boundary $\pa\O^{n-1}$ (that can be empty). Let $D^n$ be the (unbounded-)cone through $\O^{n-1}$: the region of ${\mathbb R}^n$ comprising of all vectors $\bm{x}$ whose unit radial vectors $\bm{\hat{r}}$ are in $\O^{n-1}$ (for definiteness, $r=0$ is included). Its boundary $\pa D^n$ is similarly the cone through $\pa\O^{n-1}$. Now, consider the heat equation on the cone $D^n$,
\be
\frac{\pa f}{\pa\t} = \frac{1}{2}\del^2f, \sp \del^2 = \sum_i\frac{\pa^2}{\pa x_i^2}.
\ee
$\del^2$ is the Laplacian, $\t$ is the time variable and $f(\bm{x},\bm{x}',\t)$ is the heat kernel required to satisfy the initial condition $f(\bm{x},\bm{x}',0)=\d(\bm{x}-\bm{x}')$ and a boundary condition on $\pa D^n$. Dirichlet boundary condition requires $f$ to vanish on $\pa D^n$. Neumann boundary condition requires the normal derivative of $f$ to vanish on $\pa D^n$. Standard arguments involving the Laplacian show that the solution to the heat equation is unique under these conditions and $f(\bm{x},\bm{x}',\t)$ regular at $r=0$ going to zero sufficiently fast as $r\to\inf$.

An important property of the heat kernel is that it is multiplicative in a Cartesian product. Consider a cone $D^n$ that is a Cartesian product $D_1^{n_1}\times D_2^{n_2}$ of two cones $D_1^{n_1}$ and $D_2^{n_2}$ where $n=n_1+n_2$, formed by pairing the coordinate vectors as $\bm{x}=(\bm{x}_1,\bm{x}_2)$. Then it follows from the heat equation that the heat kernel $f$ on $D^n$ is given by
\be
f(\bm{x},\bm{x}',\t) = f_1(\bm{x}_1,\bm{x}_1',\t)f_2(\bm{x}_2,\bm{x}_2',\t),
\ee
where $f_1$ and $f_2$ are the heat kernels on $D_1^{n_1}$ and $D_2^{n_2}$ respectively. The initial condition at $\t=0$ is satisfied. As can be verified, Dirichlet or Neumann boundary conditions on $D_1^{n_1}$ and $D_2^{n_2}$ imply analogous or mixed ones on $\pa D^n$ (that consists of $\pa D_1^{n_1}\times D_2^{n_2}$ and $D_1^{n_1}\times\pa D_2^{n_2}$ joined along $\pa D_1^{n_1}\times\pa D_2^{n_2}$).

There are some special cases when the heat kernel can be explicitly written down. When the cone is ${\mathbb R}$, that is the whole real line, the solution is
\be
f(x,x',\t) = \frac{1}{\sqrt{2\pi\t}}e^{-\frac{1}{2\t}\la(x-x'\ra)^2},
\lb{fs1}
\ee
that tends to $\d(x-x')$ when $\t\to 0$ as required. When the cone is ${\mathbb R}_{\ge0}$, that is the nonnegative real line, with the origin $x=0$ considered as its boundary, the solution is obtained by the method of images,
\be
f(x,x',\t) = \frac{1}{\sqrt{2\pi\t}}\la(e^{-\frac{1}{2\t}\la(x-x'\ra)^2} - e^{-\frac{1}{2\t}\la(x+x'\ra)^2}\ra) = \sqrt{\frac{2}{\pi\t}}e^{-\frac{1}{2\t}\la(x^2+x'^2\ra)}\sinh\la(\frac{xx'}{\t}\ra).
\lb{snh}
\ee
This is for Dirichlet boundary condition. For Neumann boundary condition, we will have $\cosh$ in place of $\sinh$ above.

These solutions can be combined to obtain solutions on their Cartesian product spaces making use of the multiplicative property of the heat kernel mentioned above. We thus have the heat kernel on ${\mathbb R}^n$ as the product of $n$ copies of (\ref{fs1}),
\be
f(\bm{x},\bm{x}',\t) = \frac{1}{(2\pi\t)^{\frac{n}{2}}}e^{-\frac{1}{2\t}\la(\bm{x}-\bm{x}'\ra)^2}.
\lb{fsn}
\ee
We also have the heat kernel on ${\mathbb R}_{\ge0}^n$ as the product of $n$ copies of (\ref{snh}) (under Dirichlet boundary condition),
\be
f(\bm{x},\bm{x}',\t) = \la(\frac{2}{\pi\t}\ra)^{\frac{n}{2}}e^{-\frac{1}{2\t}\la(r^2+r'^2\ra)}\prod_{i=1}^n\sinh\la(\frac{x_ix'_i}{\t}\ra).
\lb{f0xn}
\ee
We can also combine the solutions, say $p$ copies of (\ref{fs1}) and $q$ copies of (\ref{snh}), to obtain heat kernels on product spaces such as ${\mathbb R}^p\times {\mathbb R}_{\ge0}^q$.

However, the heat equation on the cone is not solvable in general, except as a series solution. One of the ways of approaching the series solution is to represent it as
\be
f(\bm{x},\bm{x}',\t) = \sum_{\n\s}\int_0^{\inf}d\m\m\,e^{-\frac{1}{2}\m^2\t}\f_{\m\n\s}(\bm{x})\f_{\m\n\s}(\bm{x}').
\lb{fff}
\ee
Here $\f_{\m\n\s}$'s are the eigenfunctions of the Laplacian on the cone $D^n$ of eigenvalue $-\m^2$, having additional labels $\n\s$ for later use,
\be
\del^2\f_{\m\n\s} = -\m^2\f_{\m\n\s}.
\lb{dl2}
\ee
It can be verified that representation (\ref{fff}) does solve the heat equation. $\f_{\m\n\s}$'s are expected to behave similarly to the heat kernel on $\pa D^n$, at $r=0$ and as $r\to\inf$. Standard arguments show that we then have $\m^2\ge0$ and that $\m^2$ covers all of the positive real line because of the unboundedness of the cone. Initial condition will be satisfied with the $\f_{\m\n\s}$'s forming a complete system, normalized according to
\be
\int_{D^n}d^nx\,\f_{\m\n\s}(\bm{x})\f_{\m'\n'\s'}(\bm{x}) = \frac{1}{\m}\d(\m-\m')\d_{\n\n'}\d_{\s\s'}.
\lb{iff}
\ee
We thus have the Laplacian boundary value problem (\ref{dl2}) to be solved that can be rewritten in radial and angular coordinates involving the spherical Laplacian.

\section{Spherical Laplacian}
\lb{Slp}

Going to radial and angular coordinates, one obtains the following expression for the Laplacian on the cone $D^n$ acting on some function $\f$,
\be
\del^2\f = \frac{1}{r^{n-1}}\frac{d}{dr}\la(r^{n-1}\frac{d\f}{dr}\ra) + \frac{1}{r^2}\del_S^2\f,
\lb{d2s}
\ee
where $\del_S^2$ acts on the angular dependence. Identifying the angular coordinates as the coordinates on the sphere $S^{n-1}$, one notes that $\del_S^2$ is the Laplace-Beltrami operator, or simply the Laplacian, on the spherical domain $\O^{n-1}$ of interest.

Angular dependence can be expressed as a dependence on $\bm{\hat{r}}=\frac{\bm{x}}{r}$ so that functions on $\O^{n-1}$ can be viewed as zero-degree (positive-)homogeneous functions of $\bm{x}$ in the cone $D^n$ (a scaling of $\bm{x}$ by say $a$ scales $r$ by $\abs{a}$ so that functions of $\bm{\hat{r}}$ remain unaffected if $a>0$). $\del_S^2$ can hence be identified with $r^2\del^2$ acting on zero-degree homogeneous functions of $\bm{x}$. Of interest to us is the Laplacian eigenvalue problem on $\O^{n-1}$,
\be
\del_{\ms S}^2h_{\n\s}(\bm{\hat{r}}) = -\l h_{\n\s}(\bm{\hat{r}}), \sp \l = \n(\n+n-2),
\lb{Lsd}
\ee
where the zero-degree homogeneous function $h_{\n\s}(\bm{\hat{r}})$, $\s$ labeling any multiplicity, is the eigenfunction satisfying appropriate boundary condition on $\pa\O^{n-1}$. As can be verified, this eigenvalue equation can be equivalently expressed as the Laplace equation $\del^2(r^{\n}h_{\n\s})=0$ in $D^n$ for a $\n$-degree homogeneous function $r^{\n}h_{\n\s}(\bm{\hat{r}})$. It is thus natural to work with the degree variable $\n$ in place of the eigenvalue $\l$.

Boundary value problems of the above kind have been extensively studied and it turns out that the eigenvalues are, and hence the $\n$'s can be taken to be, real, non-negative, unbounded and discrete, and that the eigenfunctions can be taken to be real and form a complete system. We will assume that the eigenfunctions are normalized to form an orthonormal system
\be
\int_{\O^{n-1}}d^{n-1}\hat{r}h_{\n\s}(\bm{\hat{r}})h_{\n'\s'}(\bm{\hat{r}}) = \d_{\n\n'}\d_{\s\s'},
\ee
where $d^{n-1}\hat{r}$ denotes the volume element on the unit sphere $S^{n-1}$ (area element if $S^2$).

\section{Series Solution on the Cone}
\lb{Ssc}

Going to radial and angular coordinates involving the spherical Laplacian as in (\ref{d2s}), the Laplacian boundary value problem (\ref{dl2}) on the cone $D^n$ can be rewritten as
\be
\frac{1}{r^{n-1}}\frac{d}{dr}\la(r^{n-1}\frac{d\f_{\m\n\s}}{dr}\ra) + \frac{1}{r^2}\del_S^2\f_{\m\n\s} + \m^2\f_{\m\n\s} = 0,
\lb{rfm}
\ee
where $\del_{\ms S}^2$ is the Laplacian on $\O^{n-1}$ discussed above. The complete set of eigenfunctions $h_{\n\s}(\bm{\hat{r}})$ of $\del_S^2$, of eigenvalues $\n(\n+n-2)$ and multiplicity label $\s$, enables us to look for $\f_{\m\n\s}(\bm{x})$ of the form
\be
\f_{\m\n\s}(\bm{x}) = r^{-\frac{n-2}{2}}g_{\m\n\s}(r)h_{\n\s}(\bm{\hat{r}}).
\ee
Boundary condition for $\f_{\m\n\s}(\bm{x})$ on $\pa D^n$ implies analogous one for $h_{\n\s}(\bm{\hat{r}})$ on $\pa\O^{n-1}$. Normalization (\ref{iff}) can be ensured with the $g_{\m\n\s}$'s normalized according to
\be
\int_0^{\inf}drr\,g_{\m\n\s}(r)g_{\m'\n\s}(r) = \frac{1}{\m}\d(\m-\m').
\ee
The action of $\del_{\ms S}^2$ on $h_{\n\s}(\bm{\hat{r}})$ is given by (\ref{Lsd}) so that (\ref{rfm}) gives rise to
\be
r^2\frac{\pa^2g_{\m\n\s}}{\pa r^2} + r\frac{\pa g_{\m\n\s}}{\pa r} + \la((\m r)^2-\la(\n+\frac{n-2}{2}\ra)^2\ra)g_{\m\n\s} = 0.
\ee
This is the Bessel's differential equation with the acceptable solution $g_{\m\n\s}(r)=J_{\n+\frac{n-2}{2}}\la(\m r\ra)$. Behavior of $g_{\m\n\s}(r)$ as $r\to 0$ and $r\to\inf$ is consistent with the expected one for $\f_{\m\n\s}(\bm{x})$. Orthonormality and completeness of $g_{\m\n\s}(r)$ follows from the properties of the Bessel function as the basis of the Hankel transform. We thus have
\be
\f_{\m\n\s}(\bm{x}) = r^{-\frac{n-2}{2}}J_{\n+\frac{n-2}{2}}\la(\m r\ra)h_{\n\s}(\bm{\hat{r}}).
\ee
Using this in the heat kernel representation (\ref{fff}), we get
\be
f(\bm{x},\bm{x}',\t) = (rr')^{-\frac{n-2}{2}}\sum_{\n}\int_0^{\inf}d\m\m\,e^{-\frac{1}{2}\m^2\t}J_{\n+\frac{n-2}{2}}(\m r)J_{\n+\frac{n-2}{2}}(\m r')\sum_{\s}h_{\n\s}(\bm{\hat{r}})h_{\n\s}(\bm{\hat{r}}').
\ee
Making use of an identity that expresses the integral over product of two Bessel functions in terms of the modified Bessel function, we obtain the series solution
\be
f(\bm{x},\bm{x}',\t) = \frac{1}{\t}(rr')^{-\frac{n-2}{2}}e^{-\frac{1}{2\t}(r^2+r'^2)}\sum_{\n}I_{\n+\frac{n-2}{2}}\la(\frac{rr'}{\t}\ra)\sum_{\s}h_{\n\s}(\bm{\hat{r}})h_{\n\s}(\bm{\hat{r}}').
\lb{fxn}
\ee
This result was obtained under different contexts by various authors. For $n=2$ it was obtained by \citet{sxx}. For $n=3$, it was considered within the context of circular cones by \citet{cjx}. For general dimensions, it has been presented by \citet{cxx}. The leading term in the series (\ref{fxn}) can be obtained by making use of the expansion for the Bessel functions,
\be
f(\bm{x},\bm{x}',\t) \sim \frac{2}{\G\la(\n_1+\frac{n}{2}\ra)(2\t)^{\frac{n}{2}}}\la(\frac{rr'}{2\t}\ra)^{\n_1}e^{-\frac{1}{2\t}\la(r^2+r'^2\ra)}h_{\n_1}(\bm{\hat{r}})h_{\n_1}(\bm{\hat{r}}'),
\lb{lts}
\ee
where $\n_1$ is the first $\n$ and $\G$ is the Gamma function.

In the special case of $n=2$ in polar coordinates $r$ and $\q$, we have the cone $0\le\q\le\vf$ for some $\vf<2\pi$ having the circular arc at $r=1$ as the spherical domain. For the eigenfunctions $h_{\n\s}(\bm{\hat{r}})$ under Dirichlet boundary condition, we have $\sqrt{\frac{2}{\vf}}\sin(\n\q)$ with multiplicity one and $\n=\frac{k\pi}{\vf},k=1,2,\cdots$ (under Neumann boundary condition, we have instead $\sqrt{\frac{2}{\vf}}\cos(\n\q)$ and also $\frac{1}{\sqrt{\vf}}$ for $k=0$). For $\vf=\pi$, $\n$'s are positive integers, and for $\vf=\frac{\pi}{2}$, they are positive even integers. In these cases, the heat kernel is known explicitly and the series solution is consistent with the known identities involving modified  Bessel functions.

\section{Product Domains}
\lb{Pdm}

As noted in section \ref{Hec}, an important property of the heat kernel is that it is multiplicative in a Cartesian product. Because the series solution to the heat equation involves the spectrum of the spherical Laplacian, the multiplicative property has implications on spherical domains. Towards this end, it is helpful to have a look into product domains on the sphere implied by the Cartesian products. These product domains are recognizable as topological `join's of the participating domains. In our case, they can be constructed as follows.

Let $\O_1^{n_1-1}$ on $S^{n_1-1}$ and $\O_2^{n_2-1}$ on $S^{n_2-1}$ be two spherical domains. Consider their cones $D_1^{n_1}$ in ${\mathbb R}^{n_1}$ and $D_2^{n_2}$ in ${\mathbb R}^{n_2}$. Let the product cone $D_1^{n_1}\times D_2^{n_2}$ in ${\mathbb R}^{n_1+n_2}$ be the Cartesian product of $D_1^{n_1}$ and $D_2^{n_2}$. The spherical domain associated with $D_1^{n_1}\times D_2^{n_2}$ then gives us the product domain $\O_1^{n_1-1}*\O_2^{n_2-1}$ on $S^{n_1+n_2-1}$ obtained by restricting the combined radial coordinate $r_1^2+r_2^2$ to unity. $\O_1^{n_1-1}*\O_2^{n_2-1}$ is nontrivially related to $\O_1^{n_1-1}$ and $\O_2^{n_2-1}$. Its dimension $n_1+n_2-1$ is one more than the sum of the dimensions of the two participating domains. Its boundary consists of $\pa\O_1^{n_1-1}*\O_2^{n_2-1}$ and $\O_1^{n_1-1}*\pa\O_2^{n_2-1}$ (one or both of which could be absent) joined along their boundary $\pa\O_1^{n_1-1}*\pa\O_2^{n_2-1}$.

Let $S^0$ denote the `zero-dimensional sphere': the two-point set $\{+1,-1\}$ having its `cone' as $\mathbb{R}$, the whole real line. Also, let $T^0$ be the one-point set $\{+1\}$ having its cone as $\mathbb{R}_{\ge0}$, the nonnegative real line. Then, as can be verified, $S^0*S^0$ abbreviated as $S^{0*2}$ is simply the unit circle $S^1$. More generally, $S^{n-1}$ can be obtained as $S^{0*n}$ by taking the $*$-product of $n$ copies of $S^0$. Similarly, it can be verified that $T^1\equiv T^{0*2}$ is the quadrant arc: the $90^{\circ}$-arc of the unit circle occupying one-fourth of the circumference. $T^2\equiv T^{0*3}$ is the octant triangle: the spherical triangle on the two-sphere having three $90^{\circ}$ angles occupying one-eighth of the spherical surface. $T^3\equiv T^{0*4}$ is the spherical tetrahedron on $S^3$ occupying one-sixteenth of the sphere. More generally, $T^{n-1}\equiv T^{0*n}$ is a spherical polytope on $S^{n-1}$ occupying $(2^{-n})^{\rm th}$ of the sphere. We can also construct $*-$products comprising of both $S^0$ and $T^0$. For instance, $S^0*T^0$ is the half-circle, $S^1*T^0$ is the half-sphere or half-$S^2$, and in general $S^{n-2}*T^0$ is half-$S^{n-1}$. More generally, $S^{p-1}*T^{q-1}$ is a domain on $S^{p+q-1}$ obtained by cutting away $S^{p+q-1}$ perpendicularly into half $q$-times, for instance following the steps $S^{p+q-1}\to S^{p+q-2}*T^0\to S^{p+q-3}*T^1\to\cdots\to S^{p-1}*T^{q-1}$.

Boundaries of the product domains can similarly be constructed. $\pa S^{n-1}$ is of course `zero'. Also, formally, we have $\O*\pa T^0\sim\O$ for any $\O$. This follows from the requirement that $\O*\pa T^0$, which is one of the boundary segments of $\O*T^0$ (the other being $\pa\O*T^0$), agree with $\O$ itself. Thus, $\pa T^{n-1}=\pa T^{0*n}\sim nT^{n-2}$, that is $n$ copies of $T^{n-2}$ joined at their boundaries. In the same language, we have $\pa\la(S^{p-1}*T^{q-1}\ra)\sim qS^{p-1}*T^{q-2}$.

Spherical domain $T^{n-1}$, that can be viewed as constructed out of $n$ copies of $T^0$ as $T^{0*n}$, has $n$ boundary segments of type $T^{n-2}$ that are perpendicular to each other, that is, have $90^{\circ}$ angles between their normals. It can be generalized to $T_{\bm{\r}}^{n-1}$, a similar but a `distorted' domain, having angles $\cos^{-1}\r_{ij}$ between the normals to their boundary segments $i$ and $j$. A specific case is a regular $T_{\bm{\r}}^{n-1}$, denoted $T_{(\r)}^{n-1}$, that has all those angles identical to say $\cos^{-1}\r$. As long as the $n$ normals remain linearly independent, the $n\times n$ matrix $\bm{\r}$ of $\r_{ij}$'s completed with unity along the diagonal will be positive-definite. Writing $\bm{\r}$ as $AA^T$, $A_{ik}$ being the $k^{\rm th}$ component of the $i^{\rm th}$ normal in the $\bm{x}$-coordinate system, and going to coordinates $\bm{y}=A\bm{x}$, one obtains a simpler description of the $T_{\bm{\r}}^{n-1}$ domain in terms of that of $T^{n-1}$, as one having the conical region ${\mathbb R}_{\ge0}^n$ with the boundaries determined by $y_i=0$ for the $i^{\rm th}$ boundary segment. Inner products in the $\bm{x}$-system, such as $\bm{x}^T\bm{x}'$ of two vectors $\bm{x}$ and $\bm{x}'$, can be computed in the $\bm{y}$-system as $\bm{y}^T\bm{\r}^{-1}\bm{y}'$.

\section{Spectrum On The Sphere}
\lb{SS}

Many results are known in general about the spectrum of the Laplacian on compact domains. As already stated, all eigenvalues are real, non-negative, unbounded and discrete. The first eigenvalue has multiplicity one and the corresponding eigenfunction can be taken to be positive within the domain. When the domain is the whole of the sphere $S^{n-1}$, the resulting spectrum is well-known. In this case $\n$ is an integer taking values from zero to infinity. The first $\n$, denoted $\n_1$, is zero corresponding to a constant function on $S^{n-1}$. The multiplicities of the eigenvalues will be revisited below.

On the domain $T^{n-1}$ under Dirichlet boundary condition, it is straightforward to show that $\n_1=n$. In fact, the simplest homogeneous function solving the Laplace equation in its cone ${\mathbb R}_{\ge0}^n$ and vanishing on the boundaries is of degree $n$ and is given simply by the product of the $n$-coordinates consistent with the first term in the expansion of $\sinh$ in equation (\ref{f0xn}). It is further clear that extending any spherical domain $\O^{n-1}$ by $T^0$ via a $*$-product would increase $\n_1$ by one. If extending with $S^0$ instead, $\n_1$ would remain the same. These observations are not trivial from the point of view on the sphere, but are a simple consequence of the property that $\n_1$'s are additive in a $*$-product. Additive property of $\n_1$'s in a $*$-product follows from the $\t^{-\n_1-\frac{n}{2}}$ factor in the leading term (\ref{lts}) of the series solution and the multiplicative property of the heat kernel in a Cartesian product.

To say more about the spectrum of the Laplacian on the sphere, let us next derive a spectral function, a `generating function' for the eigenvalues and multiplicities in terms of $f(\bm{x},\bm{x}',\t)$. Towards this end, let us set $\bm{x}'=\bm{x}$ and $\t=1$ to obtain
\be
f(\bm{x},\bm{x},1) = r^{2-n}e^{-r^2}\sum_{\n}I_{\n+\frac{n-2}{2}}\la(r^2\ra)\sum_{\s}(h_{\n\s}(\bm{\hat{r}}))^2.
\lb{fr2}
\ee
Note that a further operation of integrating over $\bm{x}$, along with any $\bm{\hat{r}}$-independent weight, would integrate $(h_{\n\s}(\bm{\hat{r}}))^2$ to unity (its normalization) introducing the multiplicity $m_{\n}$. This procedure derives the following expression for the spectral function $M(z)$ making use of the Laplace transform of $I_{\n}$,
\be
M(z) \equiv \sum_{\n}m_{\n}z^{\n} = \la(1-z^2\ra)z^{-\frac{n}{2}}\int_{D^n}d^nx\,e^{-\frac{1}{2z}(1-z)^2r^2}f(\bm{x},\bm{x},1),
\lb{Mzf}
\ee
where $0<z<1$ and $r$ is the length of $\bm{x}$ as given by (\ref{rxy}). If the right side can be computed, this would provide us with both the eigenvalues and the multiplicities.

The above function arose naturally from the solution of the heat equation on the cone. It differs from the usually studied trace of the heat kernel, ${\rm Tr}e^{t\del_{\ms S}^2}$, in that it is not the eigenvalues $\n(\n+n-2)$ of $-\del_{\ms S}^2$ that appear in the exponents, but rather $\n$'s themselves. Its derivation did not assume any specific character of the domain or the boundary conditions, except that $D^n$ is conical intersecting $S^{n-1}$ into some domain $\O^{n-1}$. But its applicability depends on our knowledge of $f(\bm{x},\bm{x},1)$. This is not expected to be the case in general. Below, let us first consider some special cases for which we do know $f(\bm{x},\bm{x},1)$.

Consider again for the spherical domain the whole sphere $S^{n-1}$. In this case the associated cone $D^n$ becomes the whole of ${\mathbb R}^n$. Knowing $f(\bm{x},\bm{x},1)=(2\pi)^{-\frac{n}{2}}$ from (\ref{fsn}), one readily obtains
\be
M(z) = (1-z^2)(1-z)^{-n} = 1+nz+\sum_{k=2}^{\inf}\la[{n+k-1\choose k}-{n+k-3\choose k-2}\ra]z^k.
\ee
This gives the right eigenvalues and multiplicities on the whole sphere $S^{n-1}$. The two terms inside square brackets are the dimensions of the spaces of degree $k$ and degree $k-2$ homogeneous polynomials in $n$ variables, and the role of $1-z^2$ is hence to choose the difference for the dimension of the space of degree $k$ harmonic homogeneous polynomials, that is those satisfying the Laplace equation in $n$-dimensions.

Next consider the spherical domain $T^{n-1}$. In this case, under Dirichlet boundary condition, $f(\bm{x},\bm{x},1)$ is given by (\ref{f0xn}) that gives
\be
M(z) = z^n(1-z^2)^{1-n} = \sum_{k=1}^{\inf}{n+k-3\choose k-1}z^{n+2k-2}.
\lb{mt1}
\ee
For Neumann boundary condition, factor $z^n$ in front is absent. Extended domain $T_{\bm{\r}}^{n-1}$ for $\r\ne 0$ is not solvable similarly, except when $n=2$. For $n=2$, $T_{(\r)}^1$ is a circular arc and, as we have noted earlier, $\n$'s are multiples of $\frac{\pi}{\vf}$ where $\vf=\cos^{-1}(-\r)$ and are all of multiplicity one. Its Dirichlet spectral function is hence $z^{\frac{\pi}{\vf}}(1-z^{\frac{\pi}{\vf}})^{-1}$ that becomes $z^2(1-z^2)^{-1}$ for $\r=0$, corresponding to a quadrant arc $T^1$ in agreement with (\ref{mt1}).

An important property of $M(z)$ that can be noted from its expression (\ref{Mzf}) is that, except for the factor $1-z^2$, it is multiplicative in a $*$-product as a consequence of the multiplicative property of the heat kernel in a Cartesian product. Hence, $M(z)$ on a spherical domain $\O^{n-1}$ that is a $*$-product of the form $\O_1^{n_1-1}*\O_2^{n_2-1}$ is given by
\be
M(z) = \frac{1}{1-z^2}M_1(z)M_2(z),
\lb{prd}
\ee
where $M_1(z)$ and $M_2(z)$ are the spectral functions on domains $\O_1^{n_1-1}$ and $\O_2^{n_2-1}$ respectively. It follows from above that $\n_1$'s are additive in a $*$-product (as already noted). As an example of the above product rule, consider $S^{p-1}*T^{q-1}$ for which we obtain explicitly the Dirichlet spectral function
\be
M(z) = z^{q}(1-z)^{-p}(1-z^2)^{1-q}.
\ee
For instance on the half-sphere on $S^{n-1}$ that can be constructed as $S^{n-2}*T^0$, we get $M(z)=z(1-z)^{1-n}$ which incidentally has the same multiplicities as on $T^{n-1}$. Explicit results can also be obtained with one or more $T_{(\r)}^1$'s included. Factor $z^q$ in front is absent for Neumann boundary condition. The above suggests the `atomic' results: $M(z)=1+z$ on $S^0$, $z$ on Dirichlet $T^0$ and $1$ on Neumann $T^0$.

It is in general not straightforward to obtain explicit expressions for the eigenfunctions. However, the first eigenfunction $h_{\n_1}(\bm{\hat{r}})$ (of multiplicity one) of a $*$-product domain can be obtained from those of its factor domains via the multiplicative property of the following degree-$\n_1$ homogeneous function of $\bm{x}$ in a Cartesian product ($\n_1$ being additive):
\be
\tilde{h}_{\n_1}(\bm{x}) \equiv \frac{\sqrt{2}\pi^{\frac{n}{4}}}{\sqrt{\G\la(\n_1+\frac{n}{2}\ra)}}r^{\n_1}h_{\n_1}(\bm{\hat{r}}).
\ee
This follows from the leading term of the series solution (\ref{lts}) and the multiplicative property of the heat kernel in a Cartesian product. On $S^{n-1}$ we have $\n_1=0$ and $\tilde{h}_{\n_1}(\bm{x})=1$, while on Dirichlet $T^{n-1}$ we have $\n_1=n$ and $\tilde{h}_{\n_1}(\bm{x})=2^n\prod_{i=1}^nx_i$, suggesting the `atomic' results: $\n_1=0,\tilde{h}_{\n_1}(\bm{x})=1$ on $S^0$ and $\n_1=1,\tilde{h}_{\n_1}(\bm{x})=2x$ on Dirichlet $T^0$. If interested in exploring the $h_{\n\s}(\bm{\hat{r}})$ functions in general, we could rederive our spectral function without the angular integration to obtain
\be
M(\bm{\hat{r}},\bm{\hat{r}}',z) \equiv \sum_{\n}m_{\n}(\bm{\hat{r}},\bm{\hat{r}}')z^{\n} = (1-z^2)z^{-\frac{n}{2}}\int_0^{\inf}dr\,r^{n-1}e^{-\frac{1}{2z}(1-z)^2r^2}f(r\bm{\hat{r}},r\bm{\hat{r}}',1),
\lb{Mzr}
\ee
where $m_{\n}(\bm{\hat{r}},\bm{\hat{r}}')=\sum_{\s}h_{\n\s}(\bm{\hat{r}})h_{\n\s}(\bm{\hat{r}}')$. This provides us with a spectral function for the projections on to the eigenspaces. As a function of $z\bm{\hat{r}}$ with $z$ considered as a radial coordinate, it can be identified as a kernel satisfying the Laplace equation on the cone $D^n$ for $z<1$ tending to $\d(\bm{\hat{r}}-\bm{\hat{r}}')$ as $z\to 1$ under the considered boundary condition on $\pa D^n$. In the case of ${\mathbb R}^n$ having the whole sphere $S^{n-1}$ as the spherical domain, we get
\be
M(\bm{\hat{r}},\bm{\hat{r}}',z) = \frac{1}{\abs{S^{n-1}}}\frac{\!\!\!\!1-z^2}{\la(1-2z\cos\q+z^2\ra)^{\frac{n}{2}}}, \sp \abs{S^{n-1}} = \frac{2\pi^{\frac{n}{2}}}{\G\la(\frac{n}{2}\ra)},
\ee
where $\q$ is the angle between $\bm{\hat{r}}$ and $\bm{\hat{r}}'$, and $\abs{S^{n-1}}$ is the size of the sphere $S^{n-1}$ (surface area if $S^2$). This is the Poisson kernel of the $n$-dimensional unit ball at points $z\bm{\hat{r}}$ and $\bm{\hat{r}}'$ that when expanded in powers of $z$ gives rise to zonal harmonics as projections in terms of Gegenbauer (ultraspherical) polynomials. More involved expressions can be obtained in other explicit cases such as ${\mathbb R}^p\times{\mathbb R}_{\ge0}^q$.

One may also be interested in inverting (\ref{Mzf}) to obtain information about the heat kernel on the cone when the spectral function on a domain on the sphere is known. Given $M(z)$ on a domain $\O^{n-1}$, one can obtain $f_{\O^{n-1}}(t)$, where $f_{\O^{n-1}}(r^2)=r^{n-2}\int_{\O^{n-1}}d^{n-1}\hat{r}\,f(\bm{x},\bm{x},1)$, as the Laplace inverse of $\Hat{f}_{\O^{n-1}}(s)=\frac{2z^{\frac{n}{2}}}{1-z^2}M(z),z=1+s-\sqrt{s(s+2)}$, or obtain its `trace' $\int_0^{r}drrf_{\O^{n-1}}(r^2)$ as the Laplace inverse of $\frac{1}{2s}\Hat{f}_{\O^{n-1}}(s)$. For instance, on the spherical domain $T^{n-1}$, knowing $M(z)$ from (\ref{mt1}), $f_{\O^{n-1}}(r^2)$ can be obtained as the Laplace inverse of $\frac{1}{2^{n-1}}\la(\frac{1+s-\sqrt{s(s+2)}}{s(s+2)}\ra)^{\frac{n}{2}}$. The inverse is easily carried out for $n=2$ to give
\be
\int_{\O_0^1}d\hat{r}\,f(\bm{x},\bm{x},1) = \frac{1}{4} - \frac{1}{2}I_0\la(r^2\ra)e^{-r^2} + \frac{1}{4}e^{-2r^2}.
\lb{fo2}
\ee
Here $\O_0^1$ is the quadrant arc and $I_0$ is the modified Bessel function of order zero (this can also be obtained directly from the $n=2$ series solution; alternately, knowing $M(z)$, Laplace inverse can be viewed as summing up certain series of Bessel functions). More generally, one can obtain $r^{n-2}f(r\bm{\hat{r}},r\bm{\hat{r}}',1)$ inverting (\ref{Mzr}) as the Laplace inverse of $\frac{2z^{\frac{n}{2}}}{1-z^2}M(\bm{\hat{r}},\bm{\hat{r}}',z)$.

\section{Analytical Properties}
\lb{AP}

On continuing from the $z<1$ region, $M(z)$ exhibits a singularity at $z=1$. At least for the various cases considered, the singularity is a pole of order $n-1$ (the dimension of the sphere) so that we may write around $z=1$
\be
M(z) = \frac{\!\!\!\!\!\!\!\!\,c_0}{(1-z)^{n-1}} + \frac{\!\!\!\!\!\!\!\!\!\;c_1}{(1-z)^{n-2}} + \cdots.
\lb{lau}
\ee
Coefficients $c_0$ and $c_1$ can be determined,
\be
c_0 = 2\frac{\abs{\O^{n-1}}}{\abs{S^{n-1}}}, \sp c_1 = -\frac{1}{2}c_0 - \frac{1}{2}\frac{\abs{\pa\O^{n-1}}}{\abs{S^{n-2}}}.
\lb{res}
\ee
It is convenient to write $c_1=-\frac{1}{2}(1+\g)c_0$ introducing
\be
\g = -2\frac{c_1}{c_0} - 1 = \pm\frac{1}{2}\frac{\abs{S^{n-1}}}{\abs{S^{n-2}}}\frac{\abs{\pa\O^{n-1}}}{\abs{\O^{n-1}}}.
\lb{gma}
\ee
Here, $\abs{S^{n-1}}$ and $\abs{S^{n-2}}$ are the sizes of $n-1$ and $n-2$ dimensional spheres of unit radii respectively. $\abs{\O^{n-1}}$ is the size of the domain $\O^{n-1}$ and $\abs{\pa\O^{n-1}}$ is that of its boundary $\pa\O^{n-1}$. Sizes of $\O^{n-1}$ and $\pa\O^{n-1}$ are measured in units set by the $n-1$ dimensional sphere $S^{n-1}$ of unit radius on which they reside. The sign of $\g$ is positive for Dirichlet and negative for Neumann boundary condition. Consistency with the product rule (\ref{prd}) implies that $\frac{1}{2}c_0$ is multiplicative and $\g$ is additive in a $*$-product.

The leading coefficient $c_0$ can be determined by letting $z\to 1$ in the expression for $M(z)$. Note that the exponential inside the integral would no longer provide the suppression as $r\to\inf$. As $r\to\inf$, $f(\bm{x},\bm{x},1)$ away from the boundary tends to a constant $(2\pi)^{-\frac{n}{2}}$ (see equation (\ref{fsn}) at $\bm{x}=\bm{x}',\t=1$). The integral is thus dominated by regions near $r=\inf$ where the angular integral contributes $\abs{\O^{n-1}}$. This gives, as $\e=1-z\to 0$,
\be
M(1-\e) \sim 2\e\int_0^{\inf}dr\,r^{n-1}e^{-\frac{1}{2}\e^2r^2}\frac{\abs{\O^{n-1}}}{(2\pi)^{\frac{n}{2}}} = 2\frac{\G\la(\frac{n}{2}\ra)}{2\pi^{\frac{n}{2}}}\frac{\abs{\O^{n-1}}}{\e^{n-1}}.
\ee
The factors in front can be identified as twice the inverse size of the sphere $S^{n-1}$.

The next coefficient $c_1$ can be determined by the method of images. To start with, note that the contribution to $M(z)$ coming from the source alone,
\be
\frac{c_0}{\!2}\frac{\!\!\!\!\!\!\!\!\!1+z}{(1-z)^{n-1}} = \frac{\!\!\!\!\!\!\!\!\,c_0}{(1-z)^{n-1}} - \frac{1}{2}\frac{\!\!\!\!\!\!\!\!\!\;c_0}{(1-z)^{n-2}},
\ee
makes an order $n-1$ contribution as well. In the method of images, the source placed within the domain induces images across the boundary that, under Dirichlet condition, cancel out the source effect on the boundary to ensure zero boundary condition. Since $f(\bm{x},\bm{x},1)$ is evaluated at the source location itself, as $\bm{x}$ is varied, the source moves and the images follow the source. As $r\to\inf$ many of the images will recede away from the source. The leading contribution comes from the image brought closest to the source by taking the source close to the boundary. Its contribution is $\sim-(2\pi)^{-\frac{n}{2}}e^{-2y_{\perp}^2}$. Here $y_{\perp}$ is the perpendicular distance of the source to the boundary so that the image to source distance is $2y_{\perp}$. The image contribution as $\e=1-z\to 0$ is
\be
-\frac{2\e}{(2\pi)^{\frac{n}{2}}}\int_0^{\inf}dr\,r^{n-2}e^{-\frac{1}{2}\e^2r^2}\int_{\pa\O\perp}dy_{\perp}\,e^{-2y_{\perp}^2} = -\frac{1}{2}\frac{\G\la(\frac{n-1}{2}\ra)}{2\pi^{\frac{n-1}{2}}}\frac{\abs{\pa\O^{n-1}}}{\e^{n-2}}.
\ee
The factors in front can be identified as half the inverse size of the sphere $S^{n-2}$. A negative sign is chosen to satisfy Dirichlet condition on the boundary. For Neumann boundary condition, the sign will be positive.

Expansion (\ref{lau}) is a result of an expansion of $f(\bm{x},\bm{x},1)$ in $r^{-1}$ in the expression (\ref{Mzf}) for $M(z)$. Since $\t^{\frac{n}{2}}f(\bm{x},\bm{x},\t)$ is function of the combination $\frac{r^2}{\t}$, an expansion of $f(\bm{x},\bm{x},1)$ in $r^{-1}$ is in fact an expansion of $f(\bm{x},\bm{x},\t)$ in $\sqrt{\t}$ at $\t=1$. This is the well-known expansion of the heat kernel (see for instance \citet{vxx}), in our case on the cone $D^n$. Because the higher order terms of this expansion bring in more powers of $r$ into the denominator inside the integral in (\ref{Mzf}), as such it can only be used upto coefficient $c_{n-1}$. If the remainder falls off faster than $r^{-n}$ as $r\to\inf$, its integral will be finite at $z=1$ because of the $r\to 0$ behavior of $f(\bm{x},\bm{x},1)$ evident from (\ref{fr2}). Also note here that the heat kernel expansion being an expansion in $r^{-1}$ does not see any terms of the type $e^{-r}$ for instance. That such terms are present can be seen by taking the example of the $n=2$ domain $T^1$ for which we know $\int_{\O_0^1}d\hat{r}\,f(\bm{x},\bm{x},1)$ from (\ref{fo2}). The first two terms on the right hand side of (\ref{fo2}) give rise to the heat kernel expansion while the last term, not visible to the heat kernel asymptotics, is required for the $r\to 0$ behavior.

Expansion (\ref{lau}) can also be obtained from the heat kernel expansion on $\O^{n-1}$ on the sphere itself. This can be done using the identity\footnote{Analogous relation can be written down for the pointwise object $M(\bm{\hat{r}},\bm{\hat{r}}',z)$. Inverse relations can be obtained by expressing them as Laplace transforms, giving rise to identities for the heat kernel such as the one involving the Jacobi $\q$-function on $S^1$.}
\be
M\la(e^{-s}\ra) = \frac{se^{\ell s}}{2\sqrt{\pi}}\int_0^{\inf}\frac{dt}{\;t^{\frac{3}{2}}}e^{-\ell^2t-\frac{s^2}{4t}}{\rm Tr}e^{t\del_{\ms S}^2} = \frac{e^{\ell s}}{\sqrt{\pi}}\int_0^{\inf}\frac{dt}{\!\!\sqrt{t}}e^{-t-\frac{\ell^2s^2}{4t}}{\rm Tr}e^{\frac{s^2}{4t}\del_{\ms S}^2},
\lb{Mhk}
\ee
where $\ell=\frac{1}{2}(n-2)$ and ${\rm Tr}e^{t\del_{\ms S}^2}$ is trace of the heat kernel on the sphere. The identity follows from the representation of the trace as $\sum_{\n}m_{\n}e^{-t\n(\n+2\ell)}$. Here it is convenient to expand $M\la(e^{-s}\ra)$ around $s=0$,
\be
M\la(e^{-s}\ra) \simeq \frac{\!b_0}{\,s^{n-1}} +  \frac{\!b_1}{\,s^{n-2}} +  \frac{\!b_2}{\,s^{n-3}} + \cdots.
\lb{Mes}
\ee
The heat kernel expansion on the sphere can be expressed as
\be
{\rm Tr}e^{t\del_{\ms S}^2} \simeq (4\pi t)^{-\frac{n-1}{2}}\la(a_0+a_1\sqrt{t}+a_2t+\cdots\ra).
\lb{ahk}
\ee
Coefficients of this expansion are known to certain order (see for instance \citet{vxx}), the first three of which are
\be
a_0 = \abs{\O^{n-1}}, \sp a_1 = \mp\frac{\sqrt{\pi}}{\;\;\,2}\abs{\pa\O^{n-1}}, \sp a_2 = \frac{1}{6}\int_{\O^{n-1}}R + \frac{1}{3}\int_{\pa\O^{n-1}}K.
\lb{acf}
\ee
Here $R=(n-1)(n-2)$ is the scalar curvature of $S^{n-1}$ and hence that of $\O^{n-1}$, and $K$ is the trace of  the extrinsic curvature of the boundary $\pa\O^{n-1}$ relative to $\O^{n-1}$. The sign of $a_1$ is negative for Dirichlet and positive for Neumann boundary condition. Use of the above expansions in the identity (\ref{Mhk}) gives us the $b$-coefficients,
\be
b_0 = c_0, \sp \frac{b_1}{b_0} = \ell-\frac{\g}{2}, \sp \frac{2b_2}{b_0} = \frac{b_1^2}{b_0^2} - \frac{\ell^2}{n-2} - \frac{\g^2}{4} + \frac{\!\!\!\!\!a_2}{(n-2)\abs{\O^{n-1}}},
\lb{bcf}
\ee
where $\ell=\frac{1}{2}(n-2)$ and, $c_0$ and $\g$ are given by (\ref{res}) and (\ref{gma}). The expression (\ref{acf}) for $a_2$ is applicable on domains with smooth boundaries (as well as on those without boundaries). On domains that have corners on their boundaries, it is applicable with additional contributions coming from the corner regions. For instance, on a domain on $S^2$, each vertex of angle $\vf$ on the boundary would contribute $\frac{1}{6}\la(\frac{\pi^2}{\!\!\vf}-\vf\ra)$ to $a_2$ (making use of a result attributed to \citet{kxx} that can be inferred in the present framework from the trace of the heat kernel on the $n=2$ cone following the discussion at the end of section \ref{SS}). On a domain on a higher dimensional sphere, this result is applicable when integrated along the corner regions of dimension $n-3$. For instance, on $T^{n-1}$, one finds that the corner contribution to $a_2$ is $\frac{1}{4}n(n-1)(n-2)\abs{T^{n-1}}$ (contribution from $K$ being zero) consistent with the expansion of the explicit result (\ref{mt1}).

The series expansion of the kind (\ref{Mes}) is useful in estimating the growth of the spectrum at large eigenvalues. This is done with the help of a counting function
\be
W(\n) = \sum_{\n'}m_{\n'}1_{\n'\le\n},
\ee
where $1_{\n'\le\n}$ is the step-function. $W(\n)$ counts the eigenvalues including multiplicity up to $\n$. Its Laplace transform is
\be
\Tilde{W}(s) = \int_0^{\inf}d\n W(\n)e^{-s\n} = \frac{1}{s}M(e^{-s}).
\ee
As we have noted, $M(e^{-s})$ is expected to have a pole of order $n-1$ at $s=0$ giving rise to an expansion of the kind (\ref{Mes}). Here the singularity should arise from the large $\n$ behavior of $W(\n)$. One finds for $W(\n)$ the large $\n$ expansion
\be
W(\n) \sim \frac{\;\;\;b_0\n^{n-1}}{(n-1)!} +  \frac{\;\;\;b_1\n^{n-2}}{(n-2)!} +  \frac{\;\;\;b_2\n^{n-3}}{(n-3)!} + \cdots.
\lb{Wey}
\ee
Expressed in terms of the eigenvalues $\l=\n(\n+n-2)\sim \n^2$ of the Laplacian on $\O^{n-1}$, this is consistent with the Weyl scaling law (true for more general domains).

As a Dirichlet series in $s$, one expects $M(e^{-s})$ defined on the positive real $s$-axis to be analytic on the half-plane ${\rm Re}(s)>0$. Its behavior for ${\rm Re}(s)\le 0$ is less clear. Result (\ref{Mzf}) indicates naively a relation $M\la(z^{-1}\ra)=-z^{n-2}M(z)$. However, this is not expected to hold as an approach to $z^{-1}$ from $z$ along the real axis encounters the singularity at $z=1$. For the cases considered earlier, one finds instead  $M\la(z^{-1}\ra)=(-1)^{n-1}z^{n-2-\g}M(z)$ as well as $M_D\la(z^{-1}\ra)=(-1)^{n-1}z^{n-2}M_N(z)$ where subscripts refer to Dirichlet and Neumann boundary conditions. Being consistent with the product formula (\ref{prd}), these will also hold for domains factorizable into such cases. They are however restrictive to hold in general, but when one does, $M(e^{-s})$ can be expected to be analytic on the half-plane ${\rm Re}(s)<0$ (with singularities along the imaginary $s$-axis).

\section{A Scaling Procedure}
\lb{SP}

It is a result that the eigenvalues of the Laplacian do not increase as the domain is enlarged. Having dimensions of inverse coordinate squared, eigenvalues can be expected to scale accordingly, though in general approximately, suggesting that we look for a scaling procedure to estimate the eigenvalues on a spherical domain. However, applying scaling to the eigenvalues itself, as is usually done, turns out to be not satisfactory. Let us hence look for a spectral function $M(z)$ on a target domain $\O^{n-1}$ of the form (under Dirichlet boundary condition)
\be
M(z) = z^{\a}M_0(z^{\b}),
\lb{Mab}
\ee
where $M_0(z)$ is the known spectral function on a reference domain $\O_0^{n-1}$. This implies that, given the eigenvalues $\l_{0k}=\n_{0k}(\n_{0k}+n-2),k=1,2,\cdots$ of the Laplacian on $\O_0^{n-1}$, the eigenvalues $\l_k=\n_k(\n_k+n-2)$ on $\O^{n-1}$ can be estimated according to
\be
\n_k = \a + \b\n_{0k}, \sp k = 1,2,\cdots.
\lb{abn}
\ee
Parameters $\a$ and $\b$ can be determined by expanding $M(z)$ and $M_0(z)$ into their series (\ref{lau}) at $z=1$ and matching the first two coefficients (\ref{res}) for the two domains,
\be
\a = \frac{1}{2}\la[\g-\b\g_0+(\b-1)(n-2)\ra], \sp \b = \la[\frac{\abs{\O_0^{n-1}}}{\abs{\O^{n-1}}}\ra]^{\frac{1}{n-1}},
\lb{spm}
\ee
where $\g$ and $\g_0$ for $\O^{n-1}$ and $\O_0^{n-1}$ are as given by (\ref{gma}). This estimation procedure can also be expressed as a scaling of the combination $\n+\frac{1}{2}(n-2-\g)$. Note that this does not change multiplicities. If $\O^{n-1}$ and $\O_0^{n-1}$ are closely related and the eigenvalues are well separated, this may be a reasonable assumption to make; at least for the first few eigenvalues. Eigenfunctions will of course be different.

The linear scaling estimate (\ref{abn}) with its two parameters is able to match the first two coefficients of the series expansion (\ref{lau}). It is helpful to have it extended to involve three parameters to be able to match the third coefficient of the series as well. A convenient approach to deriving such an extension without affecting multiplicities is to require the counting function $W(\n)$ to agree for the two comparison domains. Keeping the first two terms of the expansion (\ref{Wey}), $W(\n)$ can be approximated for large $\n$ as $\sim\frac{b_0}{(n-1)!}(\n+p)^{n-1}$ with $p=\frac{b_1}{b_0}$. To this order, one then obtains the scaling procedure (\ref{abn}) as a scaling of the linear combination $\n+p$. To include the next order term, let us look for a quadratic combination $(\n+p)^2+q$ for some $p$ and $q$ such that for large $\n$
\be
W(\n) \sim \frac{b_0}{(n-1)!}\la[(\n+p)^2+q\ra]^{\frac{n-1}{2}}.
\ee
Expanding this, comparing with (\ref{Wey}) and using the $b$-coefficients from (\ref{bcf}), one obtains
\be
p = \ell-\frac{\g}{2}, \sp q = -\ell^2 - \frac{1}{4}(n-2)\g^2 + \frac{\!\!\!\!\!a_2}{\abs{\O^{n-1}}},
\lb{pqa}
\ee
where $\ell=\frac{1}{2}(n-2)$, $\g$ is as defined in (\ref{gma}) and $a_2$ is given in (\ref{acf}). Thus, given $p$ and $q$ as above, $(\n+p)^2+q$ for a domain $\O^{n-1}$ can be estimated as $\b^2$ times the same combination for a reference domain $\O_0^{n-1}$. Scaling factor $\b$ is the same as before given by (\ref{spm}).

The quadratic scaling procedure thus defined is quite general applicable to domains with smooth boundaries under Dirichlet boundary condition. It can also be used on domains without boundaries in which case it becomes a linear procedure applied to the eigenvalues of the Laplacian itself. On domains that have corners on their boundaries, the procedure is applicable with additional contributions to $a_2$ coming from the corner regions as discussed in section \ref{AP}. Under Neumann boundary condition, the scaling procedure needs to be modified so as to preserve $\n_1=0$. The equivalent of the linear procedure in this case is to scale $\n(\n+2p)$ by $\b^2$. The equivalent of the quadratic procedure turns out to be to scale $(\n+p)^3+\frac{3}{2}q\n-p^3$ by $\b^3$. Both $p$ and $q$ as well as $\b$ are given by the same expressions as those obtained above (note however that $\g$ is now of the opposite sign).

\section{Numerical Comparisons}
\lb{NC}

The following numerical comparisons are for domains on the two-sphere of unit radius under Dirichlet boundary conditions, carried out with the linear scaling procedure (\ref{abn}). For clarity, area $\abs{\O^2}$ is denoted as $A$ and the perimeter $\abs{\pa\O^2}$ as $L$ so that for the scaling parameters (\ref{spm}), we have $\g=\frac{L}{\!A},\g_0=\frac{L_0}{A_0},\b=\sqrt{\frac{A_0}{\!\!A}},\a=\frac{1}{2}[\g-1-\b(\g_0-1)]$.

For the domain $T_{(\r)}^{n-1}$, the domain size expression is given by (\ref{ttr}). As $\r\to 1$, $T_{(\r)}^{n-1}$ tends to cover half the sphere. Since $M(z)=z(1-z)^{1-n}$ of the half-sphere is exactly related by scaling to $M_0(z)=z^n(1-z^2)^{1-n}$ of the $\r=0$ domain $T^{n-1}$, scaling estimates can be expected to be reasonable for $\r$ in-between. Domain $T_{(\r)}^2$ on the two-sphere is a spherical equilateral triangle of vertex angle $\cos^{-1}(-\r)$. It has $A=3\cos^{-1}(-\r)-\pi$ and $L=3\cos^{-1}\la(-\frac{\r}{1+\r}\ra)$. In this case, the reference domain for scaling estimation can be chosen to be the octant triangle $T^2$ having $M_0(z)=z^3(1-z^2)^{-2}$, $A_0=\frac{\pi}{2},L_0=\frac{3\pi}{2},\g_0=\n_{01}=3$. For $\r=\frac{1}{2}$, the spherical equilateral triangle is a tetrahedral triangle having $A=\pi,L=3\cos^{-1}\la(-\frac{1}{3}\ra)$. \citet{rtx}, in the context of a capture problem, present a theoretical and numerical result $5.159$ for the first eigenvalue $\l_1=\n_1(\n_1+1)$ of the Laplacian on the tetrahedral triangle. Scaling estimate gives $\n_1=1.826$ and $\l_1=5.162$ in excellent agreement with their result.

A spherical cap is a disk-like domain on the two-sphere. If its radius relative to its center in angles is $\q$, it has $A=2\pi(1-\cos\q)$ and $L=2\pi\sin\q$. In this case the reference domain can be chosen to be the half-sphere that has $M_0(z)=z(1-z)^{-2},\g_0=\n_{01}=1$ so that
\be
\nu_1 = \frac{1}{2}\la(\cot\frac{\q}{2}-1\ra) + \frac{\nu_{01}}{\sqrt{2}\sin\frac{\q}{2}}.
\lb{scp}
\ee
The usual scaling procedure applied to the eigenvalues of the Laplacian itself is based on just the size of the domain, and hence is not able to differentiate the effects of the boundary. \citet{rtx} present a theoretical result $\l_1=4.936$ for the first eigenvalue on a spherical cap ($\q=\frac{\pi}{3}$) having the same area as the tetrahedral triangle. Scaling with (\ref{scp}) gives $\l_1=4.949$ in excellent agreement.

A sector of the spherical cap making an angle $\vf$ at its center has $A=\vf(1-\cos\q)$ and $L=\vf\sin\q+2\q$. Choosing the reference domain to be such a sector on the half-sphere (a product domain $T_{(-\cos\vf)}^1*T^0$) that has $M_0(z)=z^{1+\frac{\pi}{\vf}}(1-z^2)^{-1}\la(1-z^{\frac{\pi}{\vf}}\ra)^{-1},\g_0=\n_{01}=1+\frac{\pi}{\vf}$, we get
\be
\nu_1 = \frac{1}{2}\la(\cot\frac{\q}{2} + \frac{\q}{\vf\sin^2\frac{\q}{2}}-\frac{\pi}{\sqrt{2}\vf\sin\frac{\q}{2}}-1\ra) + \frac{\nu_{01}}{\sqrt{2}\sin\frac{\q}{2}}.
\ee
\citet{rtx} present a theoretical result $\l_1=5.0046$ for the case $\vf=\frac{2\pi}{3}$ and $\q=\cos^{-1}\la(\frac{-1}{\sqrt{3}}\ra)$ whereas the scaling procedure gives $\l_1=5.1046$.

As a domain on the sphere is shrunk retaining its shape, it tends to approximate a flat domain in the limit, allowing for a comparison to the available solutions on flat domains. For instance, as the spherical cap has its radius $\q\to\d\sim 0$, its $\n_1\to (1+\sqrt{2})\d^{-1}=2.4142\d^{-1}$ that compares well with the flat disk solution $\sqrt{\l_1}=j_{0,1}\d^{-1}=2.4048\d^{-1}$ ($j_{0,1}$ being the first zero of the Bessel function $J_0$). The second one $\n_2\sim(1+2\sqrt{2})\d^{-1}=3.8284\d^{-1}$ also compares well with $\sqrt{\l_2}=j_{1,1}\d^{-1}=3.8317\d^{-1}$. The next one $\n_3\sim 5.2426$ is close to $\sqrt{\l_3}=5.1356$. As expected, higher ones start showing up significant differences.

Complete solution on the equilateral triangle on the plane was obtained by \citet{lxx}. Comparing the octant triangle on the sphere, one finds for the equilateral triangle of side length $\d$ on the plane $\n_1\sim \la(2\sqrt{3}+2\sqrt{\frac{2\pi}{\!\!\sqrt{3}}}\ra)\d^{-1}=7.273\d^{-1}$ that compares well with Lam\'{e}'s result $\sqrt{\l_1}=\frac{4\pi}{\!\!\sqrt{3}}\d^{-1}=7.255\d^{-1}$. The second one $\n_2\sim\la(2\sqrt{3}+4\sqrt{\frac{2\pi}{\!\!\sqrt{3}}}\ra)\d^{-1}=11.083\d^{-1}$ also compares well with $\sqrt{\l_2}=\frac{4\pi\sqrt{7}}{3}\d^{-1}=11.082\d^{-1}$. The next one $\n_3\sim 14.892$ is close to $\sqrt{\l_3}=14.510$. Here too, higher ones start showing up significant differences.

More generally, one can use the flat domain solution to estimate the first few eigenvalues on a similar domain on the sphere. Given $A$ and $L$ for a domain $\O^2$ on the sphere and $A_0$ and $L_0$ for a similar domain $\O_0^2$ on the plane, one finds for $\n_k,k=1,2,\cdots$ on $\O^2$,
\be
\n_k = \frac{1}{2}\la(\frac{L}{A}-\frac{\;\;L_0}{\sqrt{AA_0}}-1\ra) + \sqrt{\frac{A_0\l_{0k}}{\!\!\!A}},
\ee
where $\l_{0k},k=1,2,\cdots$ are the eigenvalues of the Laplacian on $\O_0^2$. This may be viewed as providing a curvature correction to the flat space eigenvalues. Note that the length scale on $\O_0^2$ cancels out, and that $A,L$ and $\nu_1$ are in units set by the unit sphere.

The above comparisons are done using the linear scaling procedure. Better results can be obtained with the quadratic scaling procedure discussed before. The quadratic procedure can be used on domains with smooth boundaries, for instance on a spherical cap on $S^{n-1}$ ($K=(n-2)\cot\q$, where $\q$ is its radius in angles), but it turns out to be identical to the linear one in the case of the spherical cap on $S^2$. As noted earlier, the quadratic procedure is also applicable on domains that have corners on their boundaries. For our numerical comparisons, this offers the improvements: $\l_1=5.1625\to 5.1606$ for the tetrahedral triangle, $\l_1=5.1046\to 5.0187$ for the sector of a spherical cap and $\sqrt{\l_1}=7.2734\d^{-1}\to 7.2613\d^{-1}$ for the equilateral triangle on the plane.

Viewed as an extension of the scaling of the Laplacian eigenvalues on flat domains to include curvature effects, the discrepancies in the scaling estimates could become significant as domains get too large relative to say the half-sphere. Also, the scaling procedure based on just two or three parameters is not expected to yield good results for all the eigenvalues, but its potential to do so for the first few is intriguing, especially because it is based on the first few coefficients of the series that governs the growth of the spectrum at large eigenvalues. It will be interesting to study the applicability of a similar scaling procedure to the spectrum of other differential operators, or to more general domains extended to a cone or by taking (\ref{Mhk}) as defining $M(z)$.

\appendix

\section{Spherical Domain Sizes}
\lb{Sdz}

Size of a spherical domain $\O^{n-1}$, denoted $\abs{\O^{n-1}}$, can be computed as the $t=1$ value of
\be
Z(t) = 2\int_{D^n}d^nx\,\d\la(r^2-t\ra).
\ee
A factor 2 is introduced since $\d(r^2-1)=\frac{1}{2}\d(r-1)$. The integral over $\bm{x}$ reduces to an integral over the angular variables giving rise to the size of the spherical domain at $t=1$. Laplace transform of $Z(t)$, denoted $\Tilde{Z}(s)$, is
\be
\Tilde{Z}(s) = \frac{2}{s^{\frac{n}{2}}}\int_{D^n}d^nx\,e^{-r^2}.
\ee
Conical nature of $D^n$ allowed us to scale $\bm{x}$ by $s^{-\frac{1}{2}}$. Inverting the transform, we get
\be
\abs{\O^{n-1}} = \abs{S^{n-1}}\frac{1}{\pi^{\frac{n}{2}}}\int_{D^n}d^nx\,e^{-r^2}, \sp \abs{S^{n-1}}=\frac{2\pi^{\frac{n}{2}}}{\G\la(\frac{n}{2}\ra)}.
\lb{Osz}
\ee
This gives us the size-fraction taken as a fraction of the size of the sphere of the same dimension. As can be noted from above, size-fractions are multiplicative in a $*$-product. For the spherical domain $T^{n-1}$, $D^n={\mathbb R}_{\ge0}^n$ and the above expression gives $2^{-n}$ for its size-fraction consistent with its description given earlier. This also follows from the multiplicative property of size-fractions given that, in this setup, we have $\abs{S^0}=2$ and $\abs{T^0}=1$ so that the size-fraction of $T^0$ can be taken to be $\frac{1}{2}$. Similarly, the size-fraction of $S^{p-1}*T^{q-1}$ turns out to be $2^{-q}$ as expected.

As noted in section \ref{Pdm}, spherical domain $T^{n-1}$ can be generalized to $T_{\bm{\r}}^{n-1}$ having elements of a matrix $\bm{\r}$ for the cosine of angles between the normals to their $n$ boundary segments. Writing $\bm{\r}$ as $AA^T$ and going from coordinates $\bm{x}$ to $\bm{y}=A\bm{x}$, we obtain from (\ref{Osz}) the domain size expression
\be
\abs{T_{\bm{\r}}^{n-1}} = \frac{\abs{S^{n-1}}}{\pi^{\frac{n}{2}}\sqrt{{\rm det}\bm{\r}}}\int_{{\mathbb R}_{\ge0}^n}d^ny\,e^{-\bm{y}^T\bm{\r}^{-1}\bm{y}}.
\lb{Tsz}
\ee
Region of $\bm{y}$-integration is ${\mathbb R}_{\ge0}^n$ since the boundaries are now determined by $y_i=0$ for the $i^{\rm th}$ boundary segment for each $i=1,\cdots,n$. $\abs{\pa T_{\bm{\r}}^{n-1}}$ can be computed from the above formula for each of the boundary segments, with $\bm{\r}^{-1}$ restricted to one dimension less.

An example of $T_{\bm{\r}}^{n-1}$ is a regular domain $T_{(\r)}^{n-1}$ with a single parameter $\r$ having angle $\cos^{-1}\r$ between the normals to their boundary segments. It has the $\bm{\r}$-matrix
\be
\bm{\r}_{ij} = (1-\r)\d_{ij} + \r, \sp \bm{\r}^{-1}_{ij} = \frac{1}{1-\r}\d_{ij} - \frac{\r}{(1-\r)(1+(n-1)\r)}.
\ee
This has determinant ${\rm det}\bm{\r}=(1-\r)^{n-1}(1+(n-1)\r)$. Matrix $A$ can be chosen to be
\be
A_{ij} = a\d_{ij}+b, \sp a = \sqrt{1-\r}, \sp b = \frac{1}{n}\la(\sqrt{1+(n-1)\r}-\sqrt{1-\r}\ra).
\ee
For this domain, the domain size expression (\ref{Tsz}) simplifies to ($\r\ge 0$)
\be
\abs{T_{(\r)}^{n-1}} = \abs{T^{n-1}}\frac{1}{\sqrt{\pi}}\int_{-\inf}^{\inf}du\,e^{-u^2}\la[{\rm erfc}\la(\frac{\sqrt{\r}\,u}{\sqrt{1-\r}}\ra)\ra]^n,
\lb{ttr}
\ee
where $\abs{T^{n-1}}=2^{-n}\abs{S^{n-1}}$ and ${\rm erfc}$ is the complementary error function. The same expression upon setting $n\to n-1$ and $\r\to\frac{\r}{1+\r}$ gives $\frac{1}{n}\abs{\pa T_{(\r)}^{n-1}}$. It can be evaluated for $\r=\frac{1}{2}$ for any $n$ giving $\abs{T_{(1/2)}^{n-1}}=\frac{1}{n+1}\abs{S^{n-1}}$ corresponding to a domain on $S^{n-1}$ analogous to a tetrahedral triangle on $S^2$. For general $\r$, on obtains the recursive differential equation
\be
\frac{\pa}{\pa\r}f_n(\r)=\frac{n(n-1)}{\pi\sqrt{1-\r^2}}f_{n-2}\la(\frac{\r}{1+2\r}\ra), \sp f_n(\r)=\frac{\abs{T_{(\r)}^{n-1}}}{\abs{T^{n-1}}},
\ee
with $f_0(\r)=f_1(\r)=f_n(0)=1$. Solving for $n=2$ gives $\abs{T_{(\r)}^1}=\cos^{-1}(-\r)$ as expected for a circular arc $T_{(\r)}^1$ of angle $\cos^{-1}(-\r)$, and for $n=3$ gives $\abs{T_{(\r)}^2}=3\cos^{-1}(-\r)-\pi$ as expected for a spherical equilateral triangle of vertex angle $\cos^{-1}(-\r)$ (and also having $\abs{\pa T_{(\r)}^2}=3\cos^{-1}\la(-\frac{\r}{1+\r}\ra)$). For small $\r$, $f_n(\r)\approx 1+\frac{1}{\pi}n(n-1)\r+\frac{1}{2\pi^2}n(n-1)(n-2)(n-3)\r^2$. Since $f_n(\r)$ is an increasing function of $\r$, we have $\abs{T_{(\r)}^{n-1}}\ge\abs{T^{n-1}}$. As $\r\to 1$, $f_n(\r)\to 2^{n-1}$ so that $T_{(\r)}^{n-1}$ tends to cover half the sphere.

\end{document}